\documentclass[twoside]{article}
\usepackage{amsfonts}
\usepackage{amssymb}
\evensidemargin=\oddsidemargin

\newcommand\bR{\mathbb R}
\newcommand\bC{\mathbb C}
\newcommand\bZ{\mathbb Z}

\newcommand\bQ{\mathbb Q}
\newcommand\bK{\mathbb K}
\newcommand\Knc{\mathbb K^{\mathrm{non-conn}}}

\newcommand\dqu{deformation quantization}

\newcommand\onto{\twoheadrightarrow}
\newcommand\es{exact sequence}
\newcommand\Kthy{$K$-theory}
\newcommand\bZm{\mathbb Z/(m)}

\newcommand{\invlim}{\raisebox{-5pt}{$\stackrel{\textstyle\lim}
{\scriptstyle\leftarrow}$}\,}
\newcommand\series{[[\hbar]]}
\newcommand\h{\ensuremath{\hbar}}
\newtheorem{Thm}{Theorem}
\newtheorem{Df}[Thm]{Definition}
\newenvironment{Def}{\begin{Df}\normalfont}{\end{Df}}
\newtheorem{Lem}[Thm]{Lemma}
\newtheorem{Prop}[Thm]{Proposition}
\newtheorem{Remk}[Thm]{Remark}
\newenvironment{Rmk}{\begin{Remk}\normalfont}{ $\square$\end{Remk}}
\newtheorem{Examp}[Thm]{Example}
\newenvironment{Ex}{\begin{Examp}\normalfont}{ $\square$\end{Examp}}
\newtheorem{Cor}[Thm]{Corollary}
\newenvironment{demo}[1][{}]{\par\noindent{\em Proof#1}. }{$\square$\par\indent}

\newenvironment{Ntn}{\noindent{\bf Notation.}}{\par}
\begin{document}
\title{Rigidity of \Kthy\ under \dqu}
\author{Jonathan Rosenberg\thanks{Partially supported by NSF grant \#
DMS-96-25336.}}
\date{Dedicated to Calvin C.\ Moore on his 60th birthday}
\maketitle
\markboth{J.\ Rosenberg}{Rigidity of \Kthy\ under \dqu}
\begin{abstract}
Quantization, at least in some formulations, 
involves replacing some algebra of observables by a (more non-commutative)
deformed algebra.  In view of the fundamental role played by \Kthy\ 
in non-commutative geometry and topology, it is
of interest to ask to what extent \Kthy\ remains ``rigid"
under this process.  We show that some positive results can be obtained
using ideas of Gabber, Gillet-Thomason, and Suslin. From this
we derive that the algebraic \Kthy\ with finite coefficients
of a \dqu\ of the functions on a compact symplectic manifold, {\em forgetting 
the topology}, recovers the topological \Kthy\ of the manifold.
\end{abstract}

{\noindent\small {\em Key words:\/} \dqu, star-product,
algebraic \Kthy, \Kthy\ with finite coefficients,
power series ring.

}

\bigskip
\begin{Ntn} If $A$ is a ring, $\bK(A)$ will denote its (connective)
\Kthy\ spectrum, the spectrum associated to the infinite loop space
$K_0(A)\times BGL(A)^+$, where $BGL(A)^+$ is the result of applying
the Quillen $+$-construction to the classifying space of the
infinite general linear group over $A$. By definition, the (algebraic)
$K$-groups $K_i(A)$ of $A$ are (at least in positive degrees)
the homotopy groups of $\bK(A)$, and the $K$-groups of $A$ with finite
coefficients $\bZm$, $K_i(A;\,\bZm)$,
are defined (at least in positive degrees) to be
the homotopy groups of $S(\bZm)\wedge \bK(A)$, where $S(\bZm)$ is the $\bZm$
Moore spectrum. These come with universal coefficient short \es es
\[
0\to K_i(A)\otimes_{\bZ}\bZm \to K_i(A;\,\bZm)
\to \mathrm{Tor}_{\bZ}(K_{i-1}(A),\,\bZm)\to 0.
\]
(This is almost, but not quite, the definition of Browder in \cite{Browder};
for an explanation of the difference between the two definitions, see
\cite{Ros1}, pp.\ 285--286.)

In the one case below where confusion might be possible between algebraic and
topological $K$-groups, we denote these by $K_j^{\mathrm{alg}}$
and $K_j^{\mathrm{top}}$, respectively.
\end{Ntn}

Now we begin with a very general definition of (formal) \dqu. Intuitively,
this is a formal deformation of the multiplication on an ``algebra
of observables,'' the deformation parameter being identified with
``Planck's constant'' $\hbar$.
\begin{Def}\label{qdef}
 Let $\overline A_0=(A_0,\,\cdot\,)$  be an algebra over 
a commutative ring $k$ (with unit), where $A_0$ is the underlying $k$-module
of $\overline A_0$ and $\cdot$ is the multiplication in $A$. A (formal)
{\bf \dqu} of $\overline A_0$ will mean an (associative) 
algebra $A=(A_0\series,\,\star\,)$ over $k\series$ (the commutative
ring of formal power series over $k$ in a variable $\h$)
with underlying $k\series$-module $A_0\series$, where the multiplication 
$\star$ in $A$ is defined by perturbing the multiplication $\cdot$
in $\overline A_0$ to a new product $\star$ via
\begin{equation}\label{starproduct}
a\star b = a\cdot b +\h\phi_1(a,b)+\h^2\phi_2(a,b)+\cdots,\qquad
a,\,b\in A_0,
\end{equation}
and then extending to series in the obvious way:
$$\left(\sum_{j=0}^\infty a_j\h^j\right)\star
\left(\sum_{l=0}^\infty b_l\h^l\right)=\sum_{j,\,l=0}^\infty
\h^{j+l}\left(a_j\cdot b_l+\sum_{p=1}^\infty \h^p\phi_p(a_j, b_l)\right).$$
Here $\phi_j,\,j=1,\,2,\,\ldots$ are $k$-bilinear maps $A_0\times A_0
\to A_0$. Note that $\overline A_0\cong A/(\h)$ {\em as algebras}, so
that one has a natural algebra map $e_0:A\onto 
\overline A_0$ (``setting \h\ to $0$''). We call the map
$e_0$  the {\bf classical limit map}.
\end{Def}
\begin{Ex}\label{trivex} A trivial but still important example is the
case where the multiplication on $\overline A_0$ is undeformed. In this
case $A=\overline A_0\series$ is simply a ring of formal power series
in one variable over the ring $\overline A_0$.
\end{Ex}
\begin{Ex}\label{symplectic}
In one of the most important examples, $k=\bC$ and $\overline A_0
=C^\infty(M)$, where $M$ is a symplectic manifold. Then there exist
non-commutative algebras $A$ satisfying Definition \ref{qdef} for which
$$\phi\sb1(f,\,g)-\phi\sb1(g,\,f)=\{f,\,g\},$$
where $\{\,\,,\,\,\}$ is the Poisson bracket on $M$. This was shown in
\cite{WLecomte} and \cite{Fed}. The $\star$-product is obtained
by a patching procedure using the Weyl quantization
of $C\sp\infty(\bR\sp{2n})$. 
\end{Ex}
In this generality, it turns out that $K_0$ is preserved under \dqu.
\begin{Thm}\label{Kzerotheorem} Let $k$ be a commutative ring with unit, 
let $\overline A_0$ be an algebra (with unit) over $k$, let $A$ be a
\dqu\ of $\overline A_0$ in the sense of Definition \ref{qdef}, and let
$e_0:A\onto \overline A_0$ be the associated classical limit map.
Then the map $(e_0)_*:K_0(A)\to K_0(\overline A_0)$ induced by $e_0$ is 
an isomorphism.
\end{Thm}
We will need the following simple lemma.
\begin{Lem}\label{invertibility}
Let $k$ be a commutative ring with unit, let $A_0$ be a $k$-algebra,
and let  $A$ be a \dqu\ of $\overline A_0$ 
in the sense of Definition \ref{qdef}. Then an element $a=\sum_{j=0}^\infty
a_j\h^j$ of $A$ ($a_j\in A_0$) is invertible if and only if $e_0(a)=a_0$
is invertible in $\overline A_0$. Similarly, an element $a=\sum_{j=0}^{n-1}
a_j\h^j$ of $A/(\h^n)$ ($a_j\in A_0$) is invertible if and only if $e_0(a)=a_0$
is invertible in $\overline A_0$. 
\end{Lem}
\begin{demo}
The ``only if'' direction is trivial, and the ``if'' direction
in the case of $A/(\h^n)$ follows from the result for $A$. The
proof (in the case of $A$) for the ``if''
direction is the usual algorithm for inversion of power series. 
More specifically, suppose $a_0$ is invertible for the multiplication $\cdot$
in $\overline A_0$, and let $a=\sum_{j=0}^\infty a_j\h^j\in A$ ($a_j\in A_0$).
We can construct an inverse $b=\sum_{l=0}^\infty b_l\h^l$ for $a$ with
respect to the product $\star$ in $A$ by letting $b_0=a_0^{-1}$ (the inverse
of $a_0$ in $\overline A_0$) and then solving for the coefficients $b_l$
by iteration in the equation
\begin{equation}\label{inv}
1=a\star b=\left(\sum_{j=0}^\infty a_j\h^j\right)\star
\left(\sum_{l=0}^\infty b_l\h^l\right)=\sum_{j,\,l=0}^\infty
\h^{j+l}\left(a_j\cdot b_l+\sum_{p=1}^\infty \h^p\phi_p(a_j, b_l)\right).
\end{equation}
Equating coefficients of powers of $\h$ on the two sides of (\ref{inv}) gives
for each $q\ge 1$ an equation (in $A_0$)
\begin{equation}\label{iterate}
\sum_{j+l+p=q}\phi_p(a_j, b_l)=0.
\end{equation}
where for convenience we let $\phi_0(a_j, b_l)=a_j\cdot b_l$. To show
these equations are (uniquely) solvable, note that 
assuming we have solved for $b_0,\,\ldots,\,
b_{q-1}$, $q\ge 1$, (\ref{iterate}) reduces to
$$a_0\cdot b_q +\sum_{\stackrel{\scriptstyle j+l+p=q}
{\scriptstyle l\le q-1}}\phi_p(a_j, b_l)=0,$$
or
$$b_q =-\sum_{\stackrel{\scriptstyle j+l+p=q}
{\scriptstyle l\le q-1}}a_0^{-1}\cdot\phi_p(a_j, b_l).$$
Thus, by induction on $q$, (\ref{inv}) has a unique solution which is a 
right $\star$-inverse to $a$. Similarly, $a$ has a unique left $\star$-inverse.
By the usual argument, these must be equal, so $a$ is invertible in $A$.
\end{demo}
\begin{demo}[ of Theorem \ref{Kzerotheorem}]
For the injectivity, it
is enough to show that if $M$ and $N$ are (left) $A$-modules
with $M\oplus N=A^n$ for some $n$, and if $\overline A_0\otimes_{e_0} M$
and $\overline A_0\otimes_{e_0} N$ are free $\overline A_0$-modules, then 
$M$ and $N$ are free $A$-modules.  Since the kernel $(\h)$ of $e_0$ is contained
in the radical of $A$ (this follows immediately from Lemma \ref{invertibility}),
the proof of Theorem 1.3.11 in \cite{Ros1} applies without change.

The proof of surjectivity is based on a version of Hensel's Lemma.
Since $A=\invlim A/(\h^j)$ and we can replace $A$ by $M_n(A)$, the $n\times n$
matrices over $A$, if necessary, it is enough to show that for $j\ge 1$,
any idempotent
$\overline a$ in $A/(\h^j)$ can be lifted to an idempotent in
$A/(\h^{j+1})$. (Then an idempotent in $\overline A_0=A/(\h)$ defining an
element of $K_0(\overline A_0)$ can be lifted by induction to an idempotent in 
$A$.) Consider the \es\ of $k$-algebras
$$
0\to I\to A/(\h^{j+1}) \to A/(\h^j) \to 0,
$$
where as a vector space, $I=\h^jA_0$, but the multiplication on $I$
vanishes since $2j\ge j+1$. Lift the idempotent
$\overline a\in A/(\h^j)$ to any element $a\in A/(\h^{j+1})$. Then $a^2-a\in
I$. But $2\overline a -1$ is invertible in $A/(\h^j)$, hence $2a-1$ is
invertible in $A/(\h^{j+1})$ by Lemma \ref{invertibility}. Let $x=(2a-1)^{-1}
\star(a-a^2)\in I$, and observe that 
$$(a+x)^2-(a+x)=(a^2+2a\star x)-(a+x)=(a^2-a)+(2a-1)\star x=0,$$
so that $a+x$ is an idempotent lifting $\overline a$.
\end{demo}
The situation for higher \Kthy\ is more complicated.
But in the remarkable paper \cite{Sus}, based on \cite{Gabber} and
\cite{GilletThomason}, Suslin computed the \Kthy\ with finite
coefficients for certain discrete valuation rings, and used the results to
study the comparison map from algebraic to topological \Kthy\ in the case
of $\bR$ and $\bC$. Some of the same
techniques can be used to prove rigidity of algebraic \Kthy\ with finite
coefficients under \dqu. Our results are basically the same as Theorem 1
in \cite{Gabber}, but without requiring the rings involved to be
commutative.
\begin{Thm}\label{maintheorem} Let $k$ be a field of characteristic
zero, let $\overline A_0$ be an algebra (with unit) over $k$, let $A$ be a
\dqu\ of $\overline A_0$ in the sense of Definition \ref{qdef}, and let
$e_0:A\onto \overline A_0$ be the associated classical limit map.
Then $e_0$ induces isomorphisms $K_j(A;\,\bZm)\stackrel{\cong}{\to}
K_j(\overline A_0;\,\bZm)$ on \Kthy\ with finite coefficients
for any $m>1$, $j>0$.
\end{Thm}
The motto of the theorem is: {\em passage to the classical limit 
preserves \Kthy\ with finite coefficients}.
But perhaps a few words of explanation for the peculiar formulation 
are in order.
\begin{enumerate}
\item We certainly cannot expect $e_0$ to induce isomorphisms of
$K$-groups {\em integrally}, since this is false even in the  case
of Example \ref{trivex}. If $\overline A_0
=k$, $\star=\cdot$, and $A=k\series$, then $A$ is a commutative local
ring and thus (see for instance \cite{Ros1}, Corollary 2.2.6) $K_1(A)=A^\times$,
which is vastly bigger than $K_1(\overline A_0)=k^\times$, and in fact 
the kernel of the map induced by $(e_0)_*$ on $\pi_1$ may be identified
with a $k$-vector space of uncountable dimension.
\item There is some subtlety in the result since $A$ is as a $k$-vector space
an infinite product of copies of $A_0$, but the \Kthy\ groups of 
an infinite product of rings are in general {\em not\/} the products of
the $K$-groups of the factors.  For a simple counterexample, let 
$R_j=C(S^{2j})$ (the continuous complex-valued functions on a sphere), 
$j=1,\,2,\,\ldots$. By Bott periodicity, $\widetilde 
K_0(R_j)\cong \bZ$. Let $b_j\in K_0(R_j)$ have non-trivial projection
into $\widetilde K_0(R_j)$. Then the element $(b_1,\,b_2,\,\ldots)$ of
$\prod_j K_0(R_j)$ does  {\em not\/} lie in the image of $K_0(\prod_j R_j)$,
since realizing $b_j$ as a formal difference of idempotent matrices 
requires matrices of increasing
size as $j\to\infty$, so that $(b_1,\,b_2,\,\ldots)$ cannot come from 
matrices of finite size over $\prod_j R_j$.  The \Kthy\ of {\em categories\/}
does commute with infinite products \cite{Carlsson}, but for quite non-trivial
reasons. However, if $\mathcal P(R)$ denotes the category of finitely generated
projective $R$-modules for a ring $R$
(the relevant category for \Kthy\ of rings), then $\mathcal 
P(\prod_j R_j)$ is not generally equivalent to $\prod_j \mathcal P(R_j)$.
\end{enumerate}
Before giving the proof, we need two preliminaries.
\begin{Lem}\label{invlim} Let $k$ be a field, let $A_0$ be a $k$-algebra,
and let  $A$ be a \dqu\ of $\overline A_0$ 
in the sense of Definition \ref{qdef}.  Then for any $n\ge 1$, the natural
maps $GL(n,\,A/(\h^{j+1}))\to GL(n,\,A/(\h^j))$ ($j=1,\,2,\,\ldots$)
are all surjective, and $GL(n,\,A)= \invlim GL(n,\,A/(\h^j))$.
\end{Lem}
\begin{demo} This follows immediately from Lemma \ref{invertibility}, 
applied not to
$A$ but to $M_n(A)$, the $n\times n$ matrices over $A$.
\end{demo}
\begin{Prop}\label{homologyiso}
Let $k$, $\overline A_0$, and $A$ be
as in Theorem \ref{maintheorem}. Then for any integers $n,\,j\ge 1$, $m>1$,
the natural map $GL(n,\,A/(\h^j))\to  GL(n,\,\overline A_0)$  
induces an isomorphism on homology with $\bZm$ coefficients.
\end{Prop}
\begin{demo} We fix $n$ and prove this by induction on $j$. The statement
is trivially true when $j=1$. So assume $j\ge 1$
and the statement is true for $j$; we'll
prove it for $j+1$. Consider the \es\ of $k$-algebras
$$
0\to I\to A/(\h^{j+1}) \to A/(\h^j) \to 0,
$$
where as a vector space, $I=\h^jA_0$, but the multiplication on $I$
vanishes since $2j\ge j+1$. By the previous lemma,
the induced map $GL(n,\,A/(\h^{j+1}))\to GL(n,\,A/(\h^j))$ 
is surjective, and the kernel $K$ consists of matrices of the form $1+x$,
$x\in M_n(I)$. Since $I^2=0$, multiplication in $K$ is given by
$(1+x)(1+y)=1+x+y$, i.e., $K\cong M_n(I)$ with its additive group structure.
Since $k$ is of characteristic zero, $K$ is therefore isomorphic to the
underlying additive group of a $\bQ$-vector space, which is uniquely
divisible. Hence $K$ is $\bZm$-acyclic, and the Hochschild-Serre
spectral sequence for 
$$1\to K\to GL(n,\,A/(\h^{j+1}))\to GL(n,\,A/(\h^j))\to 1$$ 
collapses to give $H_\bullet(GL(n,\,A/(\h^{j+1}));\,\bZm)
\cong H_\bullet(GL(n,\,A/(\h^j));\,\bZm)$. This gives the 
inductive step.
\end{demo}
\begin{demo}[ of Theorem \ref{maintheorem}]
By Lemma
\ref{invlim}, $GL(n,\,A)= \invlim GL(n,\,A/(\h^j))$ (for any $n$). Hence
the $\bZm$-homology of $GL(n,\,A)$ can be computed from that of the
$GL(n,\,A/(\h^j))$ by the Milnor $\invlim^1$ sequence. But by
Proposition \ref{homologyiso}, the maps 
$GL(n,\,A/(\h^{j+1}))\to GL(n,\,A/(\h^j))$ are all $\bZm$-homology 
isomorphisms.  Hence the inverse system $H_\bullet(GL(n,\,A/(\h^j));\,\bZm)$
(for fixed $n$) satisfies the Mittag-Leffler criterion, and
$$H_\bullet(GL(n,\,A);\,\bZm)\cong H_\bullet(GL(n,\,\overline A_0);\,\bZm).$$
Now pass the to the limit as $n\to\infty$. We deduce that the map of groups
$GL(A)\to GL(\overline A_0)$ induces a $\bZm$-homology 
isomorphism.  Applying the classifying space functor
and the Quillen $+$-construction yields that 
$BGL(A)^+\to BGL(\overline A_0)^+$ is a $\bZm$-homology equivalence (and
of course also an infinite loop map). Now the usual connective \Kthy\ spectrum
of $A$, $\bK(A)$, is just the spectrum associated to the infinite loop
structure on $K_0(A)\times BGL(A)^+$, and \Kthy\ with finite coefficients
(in positive degrees, at least) is computed by taking the homotopy groups
of $\bK(A;\,\bZm)=S(\bZm)\wedge \bK(A)$. Combining the fact that 
$BGL(A)^+\to BGL(\overline A_0)^+$ is a $\bZm$-homology equivalence with the
fact that $K_0(A)\to K_0(\overline A_0)$ is an isomorphism (Theorem 
\ref{Kzerotheorem}), we see $\bK(A;\,\bZm)\to \bK(\overline A_0;\,\bZm)$
is a homology equivalence, hence a homotopy equivalence by the Hurewicz Theorem
(which applies to {\em connective\/} spectra). (This argument bypasses
the sort of reasoning used in \cite{Sus}, Proposition 1.5, but one could
use that here instead.) So $\pi_j(\bK(A;\,\bZm))
\stackrel{\cong}{\to} \pi_j(\bK(\overline A_0;\,\bZm))$, i.e.,
$K_j(A;\,\bZm)\stackrel{\cong}{\to} K_j(\overline A_0;\,\bZm)$, for $j>0$.
\end{demo}
\begin{Cor}{\normalfont (Cf.\ \cite{Gabber}, Theorem 1, for the commutative
case.)} If $k$ is a field of characteristic zero and if
$B$ is a $k$-algebra, then for $j>0$ and any $m>1$, $K_j(B[[t]];\,\bZm)
\stackrel{\cong}{\to} K_j(B;\,\bZm)$.
\end{Cor}
\begin{demo} Apply Theorem \ref{maintheorem} to Example \ref{trivex}.
\end{demo}
\begin{Cor} Let $M$ be a compact
symplectic manifold, let $\overline A_0=C^\infty(M)$ 
with its usual Poisson structure, and let $A$ be a  \dqu\ of
$\overline A_0$. Then for $j>0$ and any $m>1$, $K_j^{\mathrm{alg}}(A;\,\bZm)
\cong K^{-j}_{\mathrm{top}}(M;\,\bZm)$, the topological \Kthy\ of $M$
with finite coefficients.
\end{Cor}
\begin{demo} We apply our results to Example \ref{symplectic}.
By Theorem \ref{maintheorem}, 
$K_0(A)\times BGL(A)^+\to K_0(\overline A_0)\times
BGL(\overline A_0)^+$ is a $\bZm$-homotopy equivalence, so for
$j>0$, 
$$
K_j^{\mathrm{alg}}(A;\,\bZm)\cong
K_j^{\mathrm{alg}}(C^\infty(M);\,\bZm).
$$
The  group on the right
is known to coincide with $K^{-j}_{\mathrm{top}}(M;\,\bZm)$ by \cite{Fischer}.
This requires comment:
Fischer's theorem is stated for the algebra of continuous functions on a 
compact space $X$, but since the proof is sheaf-theoretic, when $X$ is a 
manifold $M$, one can replace the sheaf of germs of continuous functions by
the sheaf of germs of $C^\infty$ functions, and all the arguments go through.
The essential facts needed to make everything work are:
\begin{enumerate}
\item the local ring of germs of $C^\infty$
functions at a point of a smooth manifold is Henselian;
\item for $G$ a Lie group (in particular, for $G=GL(n,\,\bC)$),
the group $C^\infty(M,\,G)$ is a ``locally convex'' topological group
in the sense of \cite{Fischer}, that is, that it is a topological group in
the $C^\infty$ topology, and that functions $M\to G$
which are close in the $C^\infty$ topology can be joined by a smooth path;
and
\item the topological \Kthy\ of $C^\infty(M)$ coincides with that
of $C(M)$ (a well-known consequence of $C^\infty$ approximation).
\end{enumerate}
\end{demo}
\begin{Rmk}
Exactly the same statement as in Theorem \ref{maintheorem} 
works when the ground ring $k$ is a field of characteristic $p$, 
except that in this case one has to assume $(m,\,p)=1$. The only difference
in the proof is that in the proof of Proposition \ref{homologyiso}, one 
should substitute  the
fact that if $(m,\,p)=1$, then a $\bZ/(p)$-vector space (regarded as a
group under addition) is $\bZm$-acyclic. In fact one can even take $k=\bZ
[{1\over m}]$ and the argument still works (see \cite{Sus}, Lemma 1.1).
\end{Rmk}
\begin{Rmk}
In fact the connective \Kthy\ spectrum is the connective 
cover of a non-connective \Kthy\ spectrum $\Knc(A)$, whose homotopy groups 
in non-negative degrees are the same as those of $K_0(A)\times BGL(A)^+$
(in other words, the Quillen $K$-groups), and whose negative homotopy
groups are the negative $K$-groups of Bass. (One of the many constructions of
this spectrum may be found in \cite{PeW1}, and a proof that it is equivalent 
to all the other standard definitions of this spectrum may be found
in \cite{PeW2},  \S\S 5--6.) An optimal statement along the
lines of Theorem \ref{maintheorem}---I am not sure whether this is
correct or not---would thus be that
$$ (e_0)_*:\Knc(A;\,\bZm)\stackrel{\cong}{\to}\Knc(\overline A_0;\,\bZm),$$
so that one gets isomorphisms similar to those of Theorem \ref{maintheorem}
for negative \Kthy\ as well,
but we have been unable to prove this. The difficulty is that the natural
way  to deloop the equivalence of Theorem \ref{maintheorem}
would be to replace $\overline A_0$ by
$\overline B_0=\overline A_0[t,\,t^{-1}]$ and define $B$ from $B_0$ by
the obvious formula derived from (\ref{starproduct}), keeping $t$ central
in $B$. The problem is that the resulting
$B$ is {\em not\/} just $A[t,\,t^{-1}]$ (which is not
$(\h)$-adically complete), but rather its $(\h)$-adic completion, and it's
not clear what effect the completion process has on $K$-groups. Other
delooping techniques run into similar problems having to do with the failure
of products and coproducts to commute. 
\end{Rmk}

\noindent{\bf Acknowledgement.} My interest in the topic of this paper
was motivated in part by a request to me from one of the
editors of {\em Mathematical Reviews\/} to write an extended review \cite{Ros2}
of two papers by Nest and Tsygan (\cite{NT1}, \cite{NT2}) on ``algebraic
index theorems.'' The more I studied the work of Nest and Tsygan, the more
I realized that it hinges on rigidity properties for cyclic homology
under \dqu. It therefore seemed natural to search for analogous
rigidity properties for \Kthy.

\pagebreak

\bigskip
\bigskip
\begin{flushleft}
\textsc{Jonathan Rosenberg}\\
Department of Mathematics\\
University of Maryland\\
College Park, MD 20742\\
email: \texttt{jmr@math.umd.edu}
\end{flushleft}
\end{document}